\documentclass{amsart}
\usepackage{amssymb}
\usepackage{amsfonts}

\newtheorem{theorem}{Theorem}
\theoremstyle{plain}

\newtheorem{corollary}{Corollary}

\newtheorem{definition}{Definition}

\newtheorem{lemma}{Lemma}

\newtheorem{remark}{Remark}

\numberwithin{equation}{section}
\input{tcilatex}

\begin{document}
\title[Some estimates for generalized commutators]{Some estimates for
generalized commutators of rough fractional maximal and integral{\ operators
on generalized weighted Morrey spaces}}
\author{F.GURBUZ}
\address{ANKARA UNIVERSITY, FACULTY OF SCIENCE, DEPARTMENT OF MATHEMATICS,
TANDO\u{G}AN 06100, ANKARA, TURKEY }
\email{feritgurbuz84@hotmail.com}
\urladdr{}
\thanks{}
\curraddr{ }
\urladdr{}
\thanks{}
\date{}
\dedicatory{}
\thanks{}

\begin{abstract}
In this paper, we establish $BMO$ estimates for generalized commutators of
rough fractional maximal and integral{\ operators on generalized weighted
Morrey spaces}, respectively.
\end{abstract}

\maketitle

\section{Introduction and main results}

The classical Morrey spaces $M_{p,\lambda }$ have been introduced by Morrey
in \cite{Morrey} to study the local behavior of solutions of second order
elliptic partial differential equations(PDEs). In recent years there has
been an explosion of interest in the study of the boundedness of operators
on Morrey-type spaces. It has been obtained that many properties of
solutions to PDEs are concerned with the boundedness of some operators on
Morrey-type spaces. In fact, better inclusion between Morrey and H\"{o}lder
spaces allows to obtain higher regularity of the solutions to different
elliptic and parabolic boundary problems. See \cite{Chiarenza, Fazio,
Palagachev, Ragusa} for details. Moreover, various Morrey spaces are defined
in the process of study. Mizuhara \cite{Mizuhara} has introduced the
generalized Morrey spaces $M_{p,\varphi }$; Komori and Shirai \cite{Komori}
has defined the weighted Morrey spaces $L_{p,\kappa }(w)$; Guliyev \cite%
{Guliyev} and Karaman \cite{Karaman} have given a concept of generalized
weighted Morrey spaces $M_{p,\varphi }\left( w\right) $ which could be
viewed as extension of both $M_{p,\varphi }$ and $L_{p,\kappa }(w)$. The
boundedness of some operators such as Hardy--Littlewood maximal operator,
the fractional integral operator as well as the fractional maximal operator
and the Calder\'{o}n--Zygmund singular integral operator on these Morrey
spaces can be seen in \cite{Guliyev, Karaman, Komori, Mizuhara}.

Let us consider the following generalized commutator of rough fractional
integral operator:%
\begin{equation*}
T_{\Omega ,\alpha }^{A}f(x)=\int \limits_{{\mathbb{R}^{n}}}\frac{\Omega (x-y)%
}{|x-y|^{n-\alpha +m-1}}R_{m}\left( A;x,y\right) f(y)dy\qquad 0<\alpha <n
\end{equation*}%
and the corresponding generalized commutator of rough fractional maximal
operator:%
\begin{equation*}
M_{\Omega ,\alpha }^{A}f(x)=\sup_{r>0}\frac{1}{r^{n-\alpha +m-1}}\dint
\limits_{|x-y|<r}\left \vert \Omega (x-y)R_{m}\left( A;x,y\right) f(y)\right
\vert dy\qquad 0<\alpha <n,
\end{equation*}%
where $\Omega \in L_{s}(S^{n-1})\left( s>1\right) $ is homogeneous of degree
zero in ${\mathbb{R}^{n}}$, $m\in 
\mathbb{N}
$, $A$ is a function defined on ${\mathbb{R}^{n}}$ and $R_{m}\left(
A;x,y\right) $ denotes the $m$-th order Taylor series remainder of $A$ at $x$
about $y$, that is, 
\begin{equation*}
R_{m}\left( A;x,y\right) =A\left( x\right) -\dsum \limits_{\left \vert
\gamma \right \vert <m}\frac{1}{\gamma !}D^{\gamma }A\left( y\right) \left(
x-y\right) ^{\gamma },
\end{equation*}%
$\gamma =\left( \gamma _{1},\cdots ,\gamma _{n}\right) $, each $\gamma
_{i}\left( i=1,\cdots ,n\right) $ is a nonnegative integer, $\left \vert
\gamma \right \vert =\dsum \limits_{i=1}^{n}\gamma _{i}$, $\gamma !=\gamma
_{1}!\cdots \gamma _{n}!$, $x^{\gamma }=x_{1}^{\gamma _{1}}\cdots
x_{n}^{\gamma _{n}}$ and $D^{\gamma }=\frac{\partial ^{\left \vert \gamma
\right \vert }}{\partial ^{\gamma _{1}}x_{1}\cdots \partial ^{\gamma
_{n}}x_{n}}$.

For $m=1$, $T_{\Omega ,\alpha }^{A}$ and $M_{\Omega ,\alpha }^{A}$ are
obviously the commutator operators,%
\begin{eqnarray*}
\left[ A,T_{\Omega ,\alpha }\right] f\left( x\right) &=&A\left( x\right)
T_{\Omega ,\alpha }f\left( x\right) -T_{\Omega ,\alpha }\left( Af\right)
\left( x\right) \\
&=&\int \limits_{{\mathbb{R}^{n}}}\frac{\Omega (x-y)}{|x-y|^{n-\alpha }}%
\left( A\left( x\right) -A\left( y\right) \right) f(y)dy
\end{eqnarray*}%
and 
\begin{eqnarray*}
\left[ A,M_{\Omega ,\alpha }\right] f\left( x\right) &=&A\left( x\right)
M_{\Omega ,\alpha }f\left( x\right) -M_{\Omega ,\alpha }\left( Af\right)
\left( x\right) \\
&=&\sup_{r>0}\frac{1}{r^{n-\alpha }}\int \limits_{|x-y|<r}\left \vert \Omega
\left( x-y\right) \right \vert \left \vert A\left( x\right) -A\left(
y\right) \right \vert \left \vert f(y)\right \vert dy,
\end{eqnarray*}%
where rough fractional integral operator $T_{\Omega ,\alpha }$ and rough
fractional maximal operator $M_{\Omega ,\alpha }$ are defined by%
\begin{equation*}
T_{\Omega ,\alpha }f(x)=\int \limits_{{\mathbb{R}^{n}}}\frac{\Omega (x-y)}{%
|x-y|^{n-\alpha }}f(y)dy\qquad 0<\alpha <n
\end{equation*}%
and 
\begin{equation*}
M_{\Omega ,\alpha }f(x)=\sup_{r>0}\frac{1}{r^{n-\alpha }}\dint%
\limits_{|x-y|<r}\left \vert \Omega (x-y)\right \vert \left \vert f(y)\right
\vert dy\qquad 0<\alpha <n,
\end{equation*}%
The weighted $\left( L_{p},L_{q}\right) $-boundedness and weak boundedness
of the operators $T_{\Omega ,\alpha }$ and $M_{\Omega ,\alpha }$ have been
given in \cite{Ding1} and \cite{Ding2}, respectively. On the other hand, if $%
m\geq 2$, then $T_{\Omega ,\alpha }^{A}$ and $M_{\Omega ,\alpha }^{A}$ are
nontrivial generalization of the above commutators, respectively. The
weighted $\left( L_{p},L_{q}\right) $-boundedness of the operators $%
T_{\Omega ,\alpha }^{A}$ and $M_{\Omega ,\alpha }^{A}$ have been given by Wu
and Yang in \cite{Wu}. In \cite{Wu}, Wu and Yang have proved the following
result.

\begin{theorem}
\label{teo1}Suppose that $0<\alpha <n$, $1<p<\frac{n}{\alpha }$, $\frac{1}{q}%
=\frac{1}{p}-\frac{\alpha }{n}$, $\Omega $ is homogeneous of degree zero
with $\Omega \in L_{s}(S^{n-1})\left( s>1\right) $. Moreover, $\left \vert
\gamma \right \vert =m-1$, $m\geq 2$ and $D^{\gamma }A\in BMO\left( {\mathbb{%
R}^{n}}\right) $. If $s^{\prime }<p$, $w\left( x\right) ^{s^{\prime }}\in
A\left( \frac{p}{s^{\prime }},\frac{q}{s^{\prime }}\right) $, then there
exists a constant $C$, independent of $A$ and $f$, such that 
\begin{equation}
\left \Vert T_{\Omega ,\alpha }^{A}f\right \Vert _{L_{q}\left( w^{q},{%
\mathbb{R}^{n}}\right) }\leq C\dsum \limits_{\left \vert \gamma \right \vert
=m-1}\left \Vert D^{\gamma }A\right \Vert _{\ast }\left \Vert f\right \Vert
_{L_{p}\left( w^{p},{\mathbb{R}^{n}}\right) },  \label{1}
\end{equation}%
\begin{equation}
\left \Vert M_{\Omega ,\alpha }^{A}f\right \Vert _{L_{q}\left( w^{q},{%
\mathbb{R}^{n}}\right) }\leq C\dsum \limits_{\left \vert \gamma \right \vert
=m-1}\left \Vert D^{\gamma }A\right \Vert _{\ast }\left \Vert f\right \Vert
_{L_{p}\left( w^{p},{\mathbb{R}^{n}}\right) }.  \label{2}
\end{equation}
\end{theorem}

Here and in the sequel, $p^{\prime }$ always denotes the conjugate index of
any $p>1$; that is, $\frac{1}{p}+\frac{1}{p^{\prime }}=1$, and $C$ stands
for a constant which is independent of the main parameters, but it may vary
from line to line.

Let $B=B(x_{0},r_{B})$ denote the ball with the center $x_{0}$ and radius $%
r_{B}$. For a given measurable set $E$, we also denote the Lebesgue measure
of $E$ by $\left \vert E\right \vert $. For any given $\Omega \subseteq {%
\mathbb{R}^{n}}$ and $0<p<\infty $, denote by $L_{p}\left( \Omega \right) $
the spaces of all functions $f$ satisfying%
\begin{equation*}
\left \Vert f\right \Vert _{L_{p}\left( \Omega \right) }=\left( \dint
\limits_{\Omega }\left \vert f\left( x\right) \right \vert ^{p}dx\right) ^{%
\frac{1}{p}}<\infty .
\end{equation*}

We recall the definition of classical Morrey spaces $M_{p,\lambda }$ as

\begin{equation*}
M_{p,\lambda }\left( {\mathbb{R}^{n}}\right) =\left \{ f:\left \Vert f\right
\Vert _{M_{p,\lambda }\left( {\mathbb{R}^{n}}\right) }=\sup \limits_{x\in {%
\mathbb{R}^{n}},r>0}\,r^{-\frac{\lambda }{p}}\, \Vert f\Vert
_{L_{p}(B(x,r))}<\infty \right \} ,
\end{equation*}%
where $f\in L_{p}^{loc}({\mathbb{R}^{n}})$, $0\leq \lambda \leq n$ and $%
1\leq p<\infty $.

Note that $M_{p,0}=L_{p}({\mathbb{R}^{n}})$ and $M_{p,n}=L_{\infty }({%
\mathbb{R}^{n}})$. If $\lambda <0$ or $\lambda >n$, then $M_{p,\lambda }={%
\Theta }$, where $\Theta $ is the set of all functions equivalent to $0$ on $%
{\mathbb{R}^{n}}$.

We also denote by $WM_{p,\lambda }\equiv WM_{p,\lambda }({\mathbb{R}^{n}})$
the weak Morrey space of all functions $f\in WL_{p}^{loc}({\mathbb{R}^{n}})$
for which 
\begin{equation*}
\left \Vert f\right \Vert _{WM_{p,\lambda }}\equiv \left \Vert f\right \Vert
_{WM_{p,\lambda }({\mathbb{R}^{n}})}=\sup_{x\in {\mathbb{R}^{n}},\;r>0}r^{-%
\frac{\lambda }{p}}\Vert f\Vert _{WL_{p}(B(x,r))}<\infty ,
\end{equation*}%
where $WL_{p}(B(x,r))$ denotes the weak $L_{p}$-space of measurable
functions $f$ for which 
\begin{equation*}
\begin{split}
\Vert f\Vert _{WL_{p}(B(x,r))}& \equiv \Vert f\chi _{_{B(x,r)}}\Vert
_{WL_{p}({\mathbb{R}^{n}})} \\
& =\sup_{t>0}t\left \vert \left \{ y\in B(x,r):\,|f(y)|>t\right \} \right
\vert ^{1/{p}} \\
& =\sup_{0<t\leq |B(x,r)|}t^{1/{p}}\left( f\chi _{_{B(x,r)}}\right) ^{\ast
}(t)<\infty ,
\end{split}%
\end{equation*}%
where $g^{\ast }$ denotes the non-increasing rearrangement of a function $g$.

Throughout the paper we assume that $x\in {\mathbb{R}^{n}}$ and $r>0$ and
also let $B(x,r)$ denotes the open ball centered at $x$ of radius $r$, $%
B^{C}(x,r)$ denotes its complement and $|B(x,r)|$ is the Lebesgue measure of
the ball $B(x,r)$ and $|B(x,r)|=v_{n}r^{n}$, where $v_{n}=|B(0,1)|$. It is
known that $M_{p,\lambda }({\mathbb{R}^{n}})$ is an expansion of $L_{p}({%
\mathbb{R}^{n}})$.

On the other hand, Mizuhara \cite{Mizuhara} has given generalized Morrey
spaces $M_{p,\varphi }$ considering $\varphi \left( r\right) $ instead of $%
r^{\lambda }$ in the above definition of the Morrey space. Later, we have
defined the generalized Morrey spaces $M_{p,\varphi }$ with normalized norm
as follows.

\begin{definition}
$\left( \text{\textbf{Generalized Morrey space}}\right) $ Let $\varphi (x,r)$
be a positive measurable function on ${\mathbb{R}^{n}}\times (0,\infty )$
and $1\leq p<\infty $. We denote by $M_{p,\varphi }\equiv M_{p,\varphi }({%
\mathbb{R}^{n}})$ the generalized Morrey space, the space of all functions $%
f\in L_{p}^{loc}({\mathbb{R}^{n}})$ with finite quasinorm 
\begin{equation*}
\Vert f\Vert _{M_{p,\varphi }}=\sup \limits_{x\in {\mathbb{R}^{n}}%
,r>0}\varphi (x,r)^{-1}\,|B(x,r)|^{-\frac{1}{p}}\, \Vert f\Vert
_{L_{p}(B(x,r))}.
\end{equation*}%
Also by $WM_{p,\varphi }\equiv WM_{p,\varphi }({\mathbb{R}^{n}})$ we denote
the weak generalized Morrey space of all functions $f\in WL_{p}^{loc}({%
\mathbb{R}^{n}})$ for which 
\begin{equation*}
\Vert f\Vert _{WM_{p,\varphi }}=\sup \limits_{x\in {\mathbb{R}^{n}}%
,r>0}\varphi (x,r)^{-1}\,|B(x,r)|^{-\frac{1}{p}}\, \Vert f\Vert
_{WL_{p}(B(x,r))}<\infty .
\end{equation*}
\end{definition}

According to this definition, we recover the Morrey space $M_{p,\lambda }$
and weak Morrey space $WM_{p,\lambda }$ under the choice $\varphi (x,r)=r^{%
\frac{\lambda -n}{p}}$: 
\begin{equation*}
M_{p,\lambda }=M_{p,\varphi }\mid _{\varphi (x,r)=r^{\frac{\lambda -n}{p}%
}},~~~~~~~~WM_{p,\lambda }=WM_{p,\varphi }\mid _{\varphi (x,r)=r^{\frac{%
\lambda -n}{p}}}.
\end{equation*}

During the last decades various classical operators, such as maximal,
singular and potential operators have been widely investigated in classical
and generalized Morrey spaces.

Komori and Shirai \cite{Komori} have introduced a version of the weighted
Morrey space $L_{p,\kappa }(w)$, which is a natural generalization of the
weighted Lebesgue space $L_{p}(w)$, and have investigated the boundedness of
classical operators in harmonic analysis.

\begin{definition}
$\left( \text{\textbf{Weighted Morrey space}}\right) $ Let $1\leq p<\infty $%
, $0<\kappa <1$ and $w$ be a weight function. We denote by $L_{p,\kappa
}(w)\equiv L_{p,\kappa }({\mathbb{R}^{n}},w)$ the weighted Morrey space of
all classes of locally integrable functions $f$ with the norm 
\begin{equation*}
\Vert f\Vert _{L_{p,\kappa }(w)}=\sup \limits_{x\in {\mathbb{R}^{n}}%
,r>0}\,w(B(x,r))^{-\frac{\kappa }{p}}\, \Vert f\Vert
_{L_{p,w}(B(x,r))}<\infty .
\end{equation*}

Furthermore, by $WL_{p,\kappa }(w)\equiv WL_{p,\kappa }({\mathbb{R}^n},w)$
we denote the weak weighted Morrey space of all classes of locally
integrable functions $f$ with the norm 
\begin{equation*}
\|f\|_{WL_{p,\kappa}(w)} = \sup \limits_{x\in{\mathbb{R}^n}, r>0} \,
w(B(x,r))^{-\frac{\kappa}{p}} \, \|f\|_{WL_{p,w}(B(x,r))} < \infty.
\end{equation*}
\end{definition}

\begin{remark}
Alternatively, we could define the weighted Morrey spaces with cubes instead
of balls. Hence we shall use these two definitions of weighted Morrey spaces
appropriate to calculation.
\end{remark}

\begin{remark}
$(1)~$ If $w\equiv 1$ and $\kappa =\lambda /n$ with $0\leq \lambda \leq n$,
then $L_{p,\lambda /n}(1)=M_{p,\lambda }({{\mathbb{R}^{n}}})$ is the
classical Morrey spaces.

$(2)~$ If $\kappa =0,$ then $L_{p,0}(w)=L_{p}(w)$ is the weighted Lebesgue
spaces.
\end{remark}

On the other hand, the generalized weighted Morrey spaces $M_{p,\varphi
}\left( w\right) $ have been introduced by Guliyev \cite{Guliyev} and
Karaman \cite{Karaman} as follows.

\begin{definition}
$\left( \text{\textbf{Generalized weighted Morrey space}}\right) $ Let $%
1\leq p<\infty $, $\varphi (x,r)$ be a positive measurable function on ${%
\mathbb{R}^{n}}\times (0,\infty )$ and $w$ be non-negative measurable
function on ${\mathbb{R}^{n}}$. We denote by $M_{p,\varphi }(w)\equiv
M_{p,\varphi }({\mathbb{R}^{n}},w)$ the generalized weighted Morrey space,
the space of all classes of functions $f\in L_{p,w}^{loc}({\mathbb{R}^{n}})$
with finite norm 
\begin{equation*}
\Vert f\Vert _{M_{p,\varphi }(w)}=\sup \limits_{x\in {\mathbb{R}^{n}}%
,r>0}\varphi (x,r)^{-1}\,w(B(x,r))^{-\frac{1}{p}}\, \Vert f\Vert
_{L_{p,w}(B(x,r))},
\end{equation*}%
where $L_{p,w}(B(x,r))$ denotes the weighted $L_{p,w}$-space of measurable
functions $f$ for which 
\begin{equation*}
\Vert f\Vert _{L_{p,w}\left( B\left( x,r\right) \right) }\equiv \Vert f\chi
_{B\left( x,r\right) }\Vert _{L_{p,w}\left( {\mathbb{R}^{n}}\right) }=\left(
\dint \limits_{B\left( x,r\right) }|f(y)|^{p}w(y)dy\right) ^{\frac{1}{p}}.
\end{equation*}%
Furthermore, by $WM_{p,\varphi }(w)\equiv WM_{p,\varphi }({\mathbb{R}^{n}}%
,w) $ we denote the weak generalized weighted Morrey space of all classes of
functions $f\in WL_{p,w}^{loc}({\mathbb{R}^{n}})$ for which 
\begin{equation*}
\Vert f\Vert _{WM_{p,\varphi }(w)}=\sup \limits_{x\in {\mathbb{R}^{n}}%
,r>0}\varphi (x,r)^{-1}\,w(B(x,r))^{-\frac{1}{p}}\, \Vert f\Vert
_{WL_{p,w}(B(x,r))}<\infty ,
\end{equation*}%
where $WL_{p,w}(B(x,r))$ denotes the weighted weak $WL_{p,w}$-space of
measurable functions $f$ for which%
\begin{equation*}
\Vert f\Vert _{WL_{p,w}\left( B\left( x,r\right) \right) }\equiv \Vert f\chi
_{B\left( x,r\right) }\Vert _{WL_{p,w}\left( {\mathbb{R}^{n}}\right) }=\sup
\limits_{t>0}tw\left( \left \{ y\in B\left( x,r\right) {:}\left \vert
f\left( y\right) \right \vert >t\right \} \right) ^{\frac{1}{p}}<\infty .
\end{equation*}
\end{definition}

\begin{remark}
$(1)~$ If $w\equiv 1$, then $M_{p,\varphi}(1)=M_{p,\varphi}$ is the
generalized Morrey space.

$(2)~$ If $\varphi (x,r)\equiv w(B(x,r))^{\frac{\kappa -1}{p}}$, $0<\kappa
<1 $, then $M_{p,\varphi }(w)=L_{p,\kappa }(w)$ is the weighted Morrey space.

$(3)~$ If $\varphi (x,r)\equiv \nu (B(x,r))^{\frac{\kappa }{p}}w(B(x,r))^{-%
\frac{1}{p}}$, $0<\kappa <1$, then $M_{p,\varphi }(w)=L_{p,\kappa }(\nu ,w)$
is the two weighted Morrey space.

$(4)~$ If $w\equiv 1$ and $\varphi (x,r)=r^{\frac{\lambda -n}{p}}$ with $%
0\leq \lambda \leq n$, then $M_{p,\varphi }(1)=M_{p,\lambda }$ is the
classical Morrey space and $WM_{p,\varphi }(1)=WM_{p,\lambda }$ is the weak
Morrey space.

$(5)~$ If $\varphi (x,r)\equiv w(B(x,r))^{-\frac{1}{p}}$, then $M_{p,\varphi
}(w)=L_{p}(w)$ is the weighted Lebesgue space.
\end{remark}

The aim of the present paper is to investigate the boundedness of
generalized commutators of rough fractional maximal and integral{\ operators
on generalized weighted Morrey spaces}, respectively. Our main results can
be formulated as follows.

\begin{theorem}
\label{teo2}Suppose that $0<\alpha <n$, $1<p<\frac{n}{\alpha }$, $\frac{1}{q}%
=\frac{1}{p}-\frac{\alpha }{n}$, $\Omega $ is homogeneous of degree zero
with $\Omega \in L_{s}(S^{n-1})\left( s>1\right) $. Moreover, let $A$ be a
function defined on ${\mathbb{R}^{n}}$, $\left \vert \gamma \right \vert
=m-1 $, $m\geq 2$ and $D^{\gamma }A\in BMO\left( {\mathbb{R}^{n}}\right) $.
If $s^{\prime }<p$, $w\left( x\right) ^{s^{\prime }}\in A\left( \frac{p}{%
s^{\prime }},\frac{q}{s^{\prime }}\right) $, then there exists a constant $C$%
, independent of $A$ and $f$, such that 
\begin{eqnarray}
\left \Vert T_{\Omega ,\alpha }^{A}f\right \Vert _{L_{q}\left(
w^{q},B(x_{0},r)\right) } &\leq &C\dsum \limits_{\left \vert \gamma \right
\vert =m-1}\left \Vert D^{\gamma }A\right \Vert _{\ast }\left( w^{q}\left(
B(x_{0},r)\right) \right) ^{\frac{1}{q}}  \notag \\
&&\times \dint \limits_{2r}^{\infty }\left( 1+\ln \frac{t}{r}\right) \left
\Vert f\right \Vert _{L_{p}\left( w^{p},B(x_{0},t)\right) }\left(
w^{q}\left( B(x_{0},t)\right) \right) ^{-\frac{1}{q}}\frac{1}{t}dt.
\label{3}
\end{eqnarray}
\end{theorem}

\begin{theorem}
\label{teo3}Let $0<\alpha <n$, $1<p<\frac{n}{\alpha }$, $\frac{1}{q}=\frac{1%
}{p}-\frac{\alpha }{n}$, $\Omega $ be homogeneous of degree zero with $%
\Omega \in L_{s}(S^{n-1})\left( s>1\right) $. Suppose that $s^{\prime }<p$, $%
w\left( x\right) ^{s^{\prime }}\in A\left( \frac{p}{s^{\prime }},\frac{q}{%
s^{\prime }}\right) $ and the pair $(\varphi _{1},\varphi _{2})$ satisfies
the condition%
\begin{equation}
\int \limits_{r}^{\infty }\left( 1+\ln \frac{t}{r}\right) \frac{\limfunc{%
essinf}\limits_{t<\tau <\infty }\varphi _{1}(x,\tau )\left( w^{p}\left(
B\left( x,\tau \right) \right) \right) ^{\frac{1}{p}}}{\left( w^{q}\left(
B\left( x,t\right) \right) \right) ^{\frac{1}{q}}}\frac{dt}{t}\leq
C_{0}\varphi _{2}(x,r),  \label{4}
\end{equation}%
where $C_{0}$ does not depend on $x$ and $r$. If $D^{\gamma }A\in BMO\left( {%
\mathbb{R}^{n}}\right) \left( \left \vert \gamma \right \vert =m-1,m\geq
2\right) $, then there is a constant $C>0$, independent of $A$ and $f$, such
that%
\begin{equation}
\left \Vert T_{\Omega ,\alpha }^{A}f\right \Vert _{M_{q,\varphi _{2}}\left(
w^{q},{\mathbb{R}^{n}}\right) }\leq C\dsum \limits_{\left \vert \gamma
\right \vert =m-1}\left \Vert D^{\gamma }A\right \Vert _{\ast }\left \Vert
f\right \Vert _{M_{p,\varphi _{1}}\left( w^{p},{\mathbb{R}^{n}}\right) },
\label{5}
\end{equation}%
\begin{equation}
\left \Vert M_{\Omega ,\alpha }^{A}f\right \Vert _{M_{q,\varphi _{2}}\left(
w^{q},{\mathbb{R}^{n}}\right) }\leq C\dsum \limits_{\left \vert \gamma
\right \vert =m-1}\left \Vert D^{\gamma }A\right \Vert _{\ast }\left \Vert
f\right \Vert _{M_{p,\varphi _{1}}\left( w^{p},{\mathbb{R}^{n}}\right) }.
\label{6}
\end{equation}
\end{theorem}

\begin{corollary}
\label{corollary1}Let $0<\alpha <n$, $1<p<\frac{n}{\alpha }$, $\frac{1}{q}=%
\frac{1}{p}-\frac{\alpha }{n}$, $\Omega $ be homogeneous of degree zero with 
$\Omega \in L_{s}(S^{n-1})\left( s>1\right) $. Let also, $s^{\prime }<p$, $%
w\left( x\right) ^{s^{\prime }}\in A\left( \frac{p}{s^{\prime }},\frac{q}{%
s^{\prime }}\right) $ and $0<\kappa <\frac{p}{q}$. If $D^{\gamma }A\in
BMO\left( {\mathbb{R}^{n}}\right) \left( \left \vert \gamma \right \vert
=m-1,m\geq 2\right) $, then there is a constant $C>0$, independent of $A$
and $f$, such that%
\begin{equation*}
\left \Vert T_{\Omega ,\alpha }^{A}f\right \Vert _{L_{q,\frac{\kappa q}{p}%
}\left( w^{q},{\mathbb{R}^{n}}\right) }\leq C\dsum \limits_{\left \vert
\gamma \right \vert =m-1}\left \Vert D^{\gamma }A\right \Vert _{\ast }\left
\Vert f\right \Vert _{L_{p,\kappa }\left( w^{p},w^{q},{\mathbb{R}^{n}}%
\right) },
\end{equation*}%
\begin{equation*}
\left \Vert M_{\Omega ,\alpha }^{A}f\right \Vert _{L_{q,\frac{\kappa q}{p}%
}\left( w^{q},{\mathbb{R}^{n}}\right) }\leq C\dsum \limits_{\left \vert
\gamma \right \vert =m-1}\left \Vert D^{\gamma }A\right \Vert _{\ast }\left
\Vert f\right \Vert _{L_{p,\kappa }\left( w^{p},w^{q},{\mathbb{R}^{n}}%
\right) }.
\end{equation*}
\end{corollary}

In the case of $w=1$ from Theorem \ref{teo3}, we get the following new
result.

\begin{corollary}
\label{corollary2}$\left( \text{see \cite{Akbulut}}\right) $ Let $0<\alpha
<n $, $1<p<\frac{n}{\alpha }$, $\frac{1}{q}=\frac{1}{p}-\frac{\alpha }{n}$, $%
\Omega $ be homogeneous of degree zero with $\Omega \in L_{s}(S^{n-1})\left(
s>1\right) $. Suppose that $s^{\prime }<p$ and the pair $(\varphi
_{1},\varphi _{2})$ satisfies the condition%
\begin{equation*}
\int \limits_{r}^{\infty }\left( 1+\ln \frac{t}{r}\right) \frac{\limfunc{%
essinf}\limits_{t<\tau <\infty }\varphi _{1}(x,\tau )\tau ^{\frac{n}{p}}}{t^{%
\frac{n}{q}}}\frac{dt}{t}\leq C_{0}\varphi _{2}(x,r),
\end{equation*}%
where $C_{0}$ does not depend on $x$ and $r$. If $D^{\gamma }A\in BMO\left( {%
\mathbb{R}^{n}}\right) \left( \left \vert \gamma \right \vert =m-1,m\geq
2\right) $, then there is a constant $C>0$, independent of $A$ and $f$, such
that%
\begin{equation*}
\left \Vert T_{\Omega ,\alpha }^{A}f\right \Vert _{M_{q,\varphi _{2}}\left( {%
\mathbb{R}^{n}}\right) }\leq C\dsum \limits_{\left \vert \gamma \right \vert
=m-1}\left \Vert D^{\gamma }A\right \Vert _{\ast }\left \Vert f\right \Vert
_{M_{p,\varphi _{1}}\left( {\mathbb{R}^{n}}\right) },
\end{equation*}%
\begin{equation*}
\left \Vert M_{\Omega ,\alpha }^{A}f\right \Vert _{M_{q,\varphi _{2}}\left( {%
\mathbb{R}^{n}}\right) }\leq C\dsum \limits_{\left \vert \gamma \right \vert
=m-1}\left \Vert D^{\gamma }A\right \Vert _{\ast }\left \Vert f\right \Vert
_{M_{p,\varphi _{1}}\left( {\mathbb{R}^{n}}\right) }.
\end{equation*}
\end{corollary}

\section{Some preliminaries and basic lemmas}

We begin with some properties of $A_{p}\left( {\mathbb{R}^{n}}\right) $
weights which play a great role in the proofs of our main results.

A weight function is a locally integrable function on ${\mathbb{R}^{n}}$
which takes values in $(0,\infty )$ almost everywhere. For a weight function 
$w$ and a measurable set $E$, we define $w(E)=\dint \limits_{E}w(x)dx$, the
Lebesgue measure of $E$ by $|E|$ and the characteristic function of $E$ by $%
\chi _{_{E}}$. Given a weight function $w$, we say that $w$ satisfies the
doubling condition if there exists a constant $D>0$ such that for any ball $%
B $, we have $w(2B)\leq Dw(B)$. When $w$ satisfies this condition, we denote 
$w\in \Delta _{2}$, for short.

If $w$ is a weight function, we denote by $L_{p}(w)\equiv L_{p}({{\mathbb{R}%
^{n}}},w)$ the weighted Lebesgue space defined by the norm 
\begin{equation*}
\Vert f\Vert _{L_{p,w}}=\left( \dint \limits_{{{\mathbb{R}^{n}}}%
}|f(x)|^{p}w(x)dx\right) ^{\frac{1}{p}}<\infty ,~~~~\text{when }1\leq
p<\infty
\end{equation*}%
and by $\Vert f\Vert _{L_{\infty ,w}}=\limfunc{esssup}\limits_{x\in {\mathbb{%
R}^{n}}}|f(x)|w(x)$ when $p=\infty $.

We denote by $WL_{p}(w)$ the weighted weak space consisting of all
measurable functions $f$ such that%
\begin{equation*}
\Vert f\Vert _{WL_{p}(w)}=\sup \limits_{t>0}tw\left( \left \{ x\in {{\mathbb{%
R}^{n}:}}\left \vert f\left( x\right) \right \vert >t\right \} \right) ^{%
\frac{1}{p}}<\infty .
\end{equation*}

We recall that a weight function $w$ is in the Muckenhoupt's class $%
A_{p}\left( {{\mathbb{R}^{n}}}\right) $, $1<p<\infty $, if 
\begin{align}
\lbrack w]_{A_{p}}& :=\sup \limits_{B}[w]_{A_{p}(B)}  \notag \\
& =\sup \limits_{B}\left( \frac{1}{|B|}\dint \limits_{B}w(x)dx\right) \left( 
\frac{1}{|B|}\dint \limits_{B}w(x)^{1-p^{\prime }}dx\right) ^{p-1}<\infty ,
\label{7}
\end{align}%
where the supremum is taken with respect to all the balls $B$ and $\frac{1}{p%
}+\frac{1}{p^{\prime }}=1$. The expression $[w]_{A_{p}}$ is called
characteristic constant of $w$. Note that, by the H\"{o}lder's inequality,
for all balls $B$ we have 
\begin{equation}
\lbrack w]_{A_{p}}^{1/p}\geq \lbrack w]_{A_{p}(B)}^{1/p}=|B|^{-1}\Vert
w\Vert _{L_{1}(B)}^{1/p}\Vert w^{-1/p}\Vert _{L_{p^{\prime }}(B)}\geq 1.
\label{8}
\end{equation}
For $p=1$, the class $A_{1}\left( {{\mathbb{R}^{n}}}\right) $ is defined by 
\begin{equation}
\frac{1}{|B|}\dint \limits_{B}w(x)dx\leq C\inf \limits_{x\in B}w\left(
x\right) ,  \label{9}
\end{equation}%
for every ball $B\subset {{\mathbb{R}^{n}}}$. Thus, we have the condition $%
Mw(x)\leq Cw(x)$ with $[w]_{A_{1}}=\sup \limits_{x\in {\mathbb{R}^{n}}}\frac{%
Mw(x)}{w(x)}$, and also for $p=\infty $ we define $A_{\infty }=\dbigcup
\limits_{1\leq p<\infty }A_{p}$, $[w]_{A_{\infty }}=\inf \limits_{1\leq
p<\infty }[w]_{A_{p}}$ and $[w]_{A_{\infty }}\leq \lbrack w]_{A_{p}}$.

One knows that $A_{p}\subset A_{q}$ if $1\leq p<q<\infty $, and that $w\in $ 
$A_{p}$ for some $1<p<q$ if $w\in A_{q}$ with $q>1$, and also $%
[w]_{A_{p}}\leq \lbrack w]_{Aq}$.

By (\ref{7}), we have 
\begin{equation}
\left( w^{-\frac{p^{\prime }}{p}}\left( B\right) \right) ^{\frac{1}{%
p^{\prime }}}=\left \Vert w^{-\frac{1}{p}}\right \Vert _{L_{p^{\prime
}}\left( B\right) }\leq C\left \vert B\right \vert w\left( B\right) ^{-\frac{%
1}{p}}  \label{10}
\end{equation}%
for $1<p<\infty $. Note that%
\begin{equation}
\left( \limfunc{essinf}\limits_{x\in E}f\left( x\right) \right) ^{-1}=%
\limfunc{esssup}\limits_{x\in E}\frac{1}{f\left( x\right) }  \label{11}
\end{equation}%
is true for any real-valued nonnegative function $f$ and is measurable on $E$
(see \cite{Wheeden-Zygmund} page 143) and (\ref{9}); we get%
\begin{eqnarray}
\left \Vert w^{-1}\right \Vert _{L_{\infty }\left( B\right) } &=&\limfunc{%
esssup}\limits_{x\in B}\frac{1}{w\left( x\right) }  \notag \\
&=&\frac{1}{\limfunc{essinf}\limits_{x\in B}w\left( x\right) }\leq C\left
\vert B\right \vert w\left( B\right) ^{-1}.  \label{12}
\end{eqnarray}

We also need another weight class $A_{p,q}$ introduced by Muckenhoupt and
Wheeden in \cite{Muckenhoupt} to study weighted boundedness of fractional
integral{\ operators.}

A weight function $w$ belongs to the Muckenhoupt-Wheeden class $A_{p,q}$ 
\cite{Muckenhoupt} for $1<p<q<\infty $ if 
\begin{align}
\lbrack w]_{A_{p,q}}& :=\sup \limits_{B}[w]_{A_{p,q}(B)}  \notag \\
& =\sup \limits_{B}\left( \frac{1}{|B|}\dint \limits_{B}w(x)^{q}dx\right) ^{%
\frac{1}{q}}\left( \frac{1}{|B|}\dint \limits_{B}w(x)^{-p^{\prime
}}dx\right) ^{\frac{1}{p^{\prime }}}<\infty ,  \label{13}
\end{align}%
where the supremum is taken with respect to all the balls $B$. Note that, by
the H\"{o}lder's inequality, for all balls $B$ we have 
\begin{equation}
\lbrack w]_{A_{p,q}}\geq \lbrack w]_{A_{p,q}(B)}=|B|^{\frac{1}{p}-\frac{1}{q}%
-1}\Vert w\Vert _{L_{q}(B)}\Vert w^{-1}\Vert _{L_{p^{\prime }}(B)}\geq 1.
\label{14}
\end{equation}%
Moreover, if $\frac{1}{q}=\frac{1}{p}-\frac{\alpha }{n}$ with $1<p<\frac{n}{%
\alpha }$ and $0<\alpha <n$, then it's easy to deduce that%
\begin{equation*}
w\left( x\right) \in A_{p,q}\iff w\left( x\right) ^{q}\in A_{\frac{q\left(
n-\alpha \right) }{n}}\iff w\left( x\right) ^{q}\in A_{1+\frac{q}{p^{\prime }%
}}.
\end{equation*}

For $p=1$, $w$ is in $A_{1,q}$ with $1<q<\infty $ if%
\begin{eqnarray*}
\lbrack w]_{A_{1,q}} &:&=\sup \limits_{B}[w]_{A_{1,q}(B)} \\
&=&\sup \limits_{B}\left( \frac{1}{|B|}\dint \limits_{B}w(x)^{q}dx\right) ^{%
\frac{1}{q}}\left( \limfunc{esssup}\limits_{x\in B}\frac{1}{w\left( x\right) 
}\right) <\infty ,
\end{eqnarray*}%
where the supremum is taken with respect to all the balls $B$. Thus, we get%
\begin{equation}
\left( \frac{1}{|B|}\dint \limits_{B}w(x)^{q}dx\right) ^{\frac{1}{q}}\leq
C\inf \limits_{x\in B}w\left( x\right) ,  \label{15}
\end{equation}%
for every ball $B\subset {{\mathbb{R}^{n}}}$.

By (\ref{13}), we have%
\begin{equation}
\left( \dint \limits_{B}w(x)^{q}dx\right) ^{\frac{1}{q}}\left( \dint
\limits_{B}w(x)^{-p^{\prime }}dx\right) ^{\frac{1}{p^{\prime }}}\leq C\left
\vert B\right \vert ^{\frac{1}{q}+\frac{1}{p^{\prime }}}.  \label{16}
\end{equation}

We summarize some properties about Muckenhoupt-Wheeden class $A_{p,q}$; see 
\cite{GarRub, Muckenhoupt}.

\begin{lemma}
\label{lemma0}Given $1\leq p\leq q<\infty $. The following statements hold.

$\left( i\right) $ $w\left( x\right) \in A_{p,q}\iff w\left( x\right)
^{q}\in A_{1+\frac{q}{p^{\prime }}};$

$\left( ii\right) $ $w\left( x\right) \in A_{p,q}\iff w\left( x\right)
^{-p^{\prime }}\in A_{1+\frac{p^{\prime }}{q}};$

$\left( iii\right) $ If $p_{1}<p_{2}$ and $q_{1}<q_{2}$, then $A\left(
p_{1},q_{1}\right) \subset A\left( p_{2},q_{2}\right) $.
\end{lemma}

Let us recall the definition and some properties of $BMO\left( {\mathbb{R}%
^{n}}\right) $. A locally integrable function $b$ is said to be in $BMO$ if 
\begin{equation*}
\Vert b\Vert _{\ast }=\sup_{x\in {\mathbb{R}^{n}},r>0}\frac{1}{|B(x,r)|}%
\dint \limits_{B(x,r)}|b(y)-b_{B(x,r)}|dy<\infty ,
\end{equation*}%
where 
\begin{equation*}
b_{B(x,r)}=\frac{1}{|B(x,r)|}\dint \limits_{B(x,r)}b(y)dy.
\end{equation*}

Define 
\begin{equation*}
BMO({\mathbb{R}^{n}})=\{b\in L_{1}^{loc}({\mathbb{R}^{n}})~:~\Vert b\Vert
_{\ast }<\infty \}.
\end{equation*}

If one regards two functions whose difference is a constant as one (modulo
constants), then the space $BMO({\mathbb{R}^{n}})$ is a Banach space with
respect to norm $\Vert \cdot \Vert _{\ast }$.

An early work about $BMO({\mathbb{R}^{n}})$ space can be attributed to John
and Nirenberg \cite{John-Nirenberg}. For $1<p<\infty $, there is a close
relation between $BMO({\mathbb{R}^{n}})$ and $A_{p}$ weights:%
\begin{equation*}
BMO({\mathbb{R}^{n}})=\left \{ \alpha \log w:w\in A_{p}\text{, }\alpha \geq
0\right \} .
\end{equation*}

\begin{lemma}
\label{lemma1}$\left( \text{John-Nirenberg inequality; see \cite%
{John-Nirenberg}}\right) $ There are constants $C_{1}$, $C_{2}>0$, such that
for all $b\in BMO({\mathbb{R}^{n}})$ and $\beta >0$ 
\begin{equation*}
\left \vert \left \{ x\in B\,:\,|b(x)-b_{B}|>\beta \right \} \right \vert
\leq C_{1}|B|e^{-C_{2}\beta /\Vert b\Vert _{\ast }},~~~\forall B\subset {%
\mathbb{R}^{n}}.
\end{equation*}
\end{lemma}

By Lemma \ref{lemma1}, it is easy to get the following.

\begin{lemma}
\label{lemma2}Let $w\in A_{\infty }$ and $b\in BMO({\mathbb{R}^{n}})$. Then
for any $p\geq 1$ we have 
\begin{equation*}
\left( \frac{1}{w(B)}\dint \limits_{B}|b(y)-b_{B}|^{p}w(y)dy\right) ^{\frac{1%
}{p}}\leq C\Vert b\Vert _{\ast }.
\end{equation*}
\end{lemma}

\begin{lemma}
\label{lemma3}$\left( \text{see \cite{LinLu}}\right) $ Let $b$ be a function
in $BMO({\mathbb{R}^{n}})$. Let also $1\leq p<\infty $, $x\in {\mathbb{R}^{n}%
}$, and $r_{1},r_{2}>0$. Then 
\begin{equation*}
\left( \frac{1}{|B(x,r_{1})|}\dint%
\limits_{B(x,r_{1})}|b(y)-b_{B(x,r_{2})}|^{p}dy\right) ^{\frac{1}{p}}\leq
C\left( 1+\left \vert \ln \frac{r_{1}}{r_{2}}\right \vert \right) \Vert
b\Vert _{\ast },
\end{equation*}%
where $C>0$ is independent of $b$, $x$, $r_{1}$ and $r_{2}$.
\end{lemma}

By Lemma \ref{lemma2} and Lemma \ref{lemma3}, it is easily to prove the
following result.

\begin{lemma}
\label{lemma4}Let $w\in A_{\infty }$ and $b\in BMO({\mathbb{R}^{n}})$. Let
also $1\leq p<\infty $, $x\in {\mathbb{R}^{n}}$, and $r_{1},r_{2}>0$. Then 
\begin{equation}
\left( \frac{1}{w(B(x,r_{1}))}\int%
\limits_{B(x,r_{1})}|b(y)-b_{B(x,r_{2})}|^{p}w(y)dy\right) ^{\frac{1}{p}%
}\leq C\left( 1+\left \vert \ln \frac{r_{1}}{r_{2}}\right \vert \right)
\Vert b\Vert _{\ast },  \label{31}
\end{equation}%
where $C>0$ is independent of $b$, $w$, $x$, $r_{1}$ and $r_{2}$.
\end{lemma}

At the end of this section, we list some known results about $R_{m}\left(
A;x,y\right) $.

\begin{lemma}
$\left( \text{see \cite{Cohen}}\right) $ Let $A$ be a function on ${\mathbb{R%
}^{n}}$ and $D^{\gamma }A\in L_{q}^{loc}({\mathbb{R}^{n}})$ for $\left \vert
\gamma \right \vert =m$ and some $q>n$. Then%
\begin{equation*}
\left \vert R_{m}\left( A;x,y\right) \right \vert \leq C\left \vert
x-y\right \vert ^{m}\dsum \limits_{\left \vert \gamma \right \vert =m}\left( 
\frac{1}{\left \vert \tilde{Q}\left( x,y\right) \right \vert }\dint \limits_{%
\tilde{Q}\left( x,y\right) }\left \vert D^{\gamma }A\left( z\right) \right
\vert ^{q}dz\right) ^{\frac{1}{q}},
\end{equation*}%
where $\tilde{Q}\left( x,y\right) $ is the cube centered at $x$ with edges
parallel to the axes and having diameter $5\sqrt{n}\left \vert
x-y\right
\vert $.
\end{lemma}

\begin{lemma}
$\left( \text{see \cite{Cohen}}\right) $ For fixed $x\in {\mathbb{R}^{n}}$,
let%
\begin{equation*}
\overline{A}\left( x\right) =A\left( x\right) -\dsum \limits_{\left \vert
\gamma \right \vert =m-1}\frac{1}{\gamma !}\left( D^{\gamma }A\right)
_{B\left( x,5\sqrt{n}\left \vert x-y\right \vert \right) }x^{\gamma }.
\end{equation*}%
Then $R_{m}\left( A;x,y\right) =R_{m}\left( \overline{A};x,y\right) $.
\end{lemma}

\begin{lemma}
\label{lemma5}$\left( \text{see \cite{Akbulut}}\right) $ Let $x\in
B(x_{0},r) $, $y\in B\left( x_{0},2^{j+1}r\right) \diagdown B\left(
x_{0},2^{j}r\right) $. Then%
\begin{equation}
\left \vert R_{m}\left( A;x,y\right) \right \vert \leq C\left \vert
x-y\right \vert ^{m-1}\left( j\dsum \limits_{\left \vert \gamma \right \vert
=m-1}\left \Vert D^{\gamma }A\right \Vert _{\ast }+\dsum \limits_{\left
\vert \gamma \right \vert =m-1}\left \vert D^{\gamma }A\left( y\right)
-\left( D^{\gamma }A\right) _{B\left( x_{0},r\right) }\right \vert \right) .
\label{34}
\end{equation}
\end{lemma}

\section{Proofs of the main results}

\textbf{Proof of Theorem \ref{teo2}.}

We write as $f=f_{1}+f_{2}$, where $f_{1}\left( y\right) =f\left( y\right)
\chi _{B\left( x_{0},2r\right) }\left( y\right) $, $\chi _{B\left(
x_{0},2r\right) }$ denotes the characteristic function of $B\left(
x_{0},2r\right) $. Then%
\begin{equation*}
\left \Vert T_{\Omega ,\alpha }^{A}f\right \Vert _{L_{q}\left( w^{q},B\left(
x_{0},r\right) \right) }\leq \left \Vert T_{\Omega ,\alpha }^{A}f_{1}\right
\Vert _{L_{q}\left( w^{q},B\left( x_{0},r\right) \right) }+\left \Vert
T_{\Omega ,\alpha }^{A}f_{2}\right \Vert _{L_{q}\left( w^{q},B\left(
x_{0},r\right) \right) }.
\end{equation*}%
Since $f_{1}\in L_{p}\left( w^{p},{\mathbb{R}^{n}}\right) $, by the
boundedness of $T_{\Omega ,\alpha }^{A}$ from $L_{p}\left( w^{p},{\mathbb{R}%
^{n}}\right) $ to $L_{q}\left( w^{q},{\mathbb{R}^{n}}\right) $ (see Theorem %
\ref{teo1}) we get 
\begin{eqnarray*}
\left \Vert T_{\Omega ,\alpha }^{A}f_{1}\right \Vert _{L_{q}\left(
w^{q},B\left( x_{0},r\right) \right) } &\leq &\left \Vert T_{\Omega ,\alpha
}^{A}f_{1}\right \Vert _{L_{q}\left( w^{q},{\mathbb{R}^{n}}\right) } \\
&\leq &C\dsum \limits_{\left \vert \gamma \right \vert =m-1}\left \Vert
D^{\gamma }A\right \Vert _{\ast }\left \Vert f_{1}\right \Vert _{L_{p}\left(
w^{p},{\mathbb{R}^{n}}\right) } \\
&=&C\dsum \limits_{\left \vert \gamma \right \vert =m-1}\left \Vert
D^{\gamma }A\right \Vert _{\ast }\left \Vert f\right \Vert _{L_{p}\left(
w^{p},B\left( x_{0},2r\right) \right) }.
\end{eqnarray*}%
Since $1<p<q$ and $\frac{s^{\prime }p}{p^{\prime }\left( p-s^{\prime
}\right) }\geq 1$, then by the H\"{o}lder's inequality%
\begin{eqnarray*}
1 &\leq &\left( \frac{1}{|B|}\dint \limits_{B}w(y)^{p}dy\right) ^{\frac{1}{p}%
}\left( \frac{1}{|B|}\dint \limits_{B}w(y)^{-p^{\prime }}dy\right) ^{\frac{1%
}{p^{\prime }}} \\
&\leq &\left( \frac{1}{|B|}\dint \limits_{B}w(y)^{q}dy\right) ^{\frac{1}{q}%
}\left( \frac{1}{|B|}\dint \limits_{B}w(y)^{-\frac{s^{\prime }p}{\left(
p-s^{\prime }\right) }}dy\right) ^{\frac{\left( p-s^{\prime }\right) }{%
s^{\prime }p}}.
\end{eqnarray*}%
This means%
\begin{equation*}
r^{\frac{n}{s^{\prime }}-\alpha }\leq \left( w^{q}\left( B\left(
x_{0},r\right) \right) \right) ^{\frac{1}{q}}\left \Vert w^{-1}\right \Vert
_{L_{\frac{s^{\prime }p}{\left( p-s^{\prime }\right) }}\left( B\left(
x_{0},r\right) \right) }.
\end{equation*}%
Thus,%
\begin{eqnarray*}
\left \Vert f\right \Vert _{L_{p}\left( w^{p},B\left( x_{0},2r\right)
\right) } &\leq &Cr^{\frac{n}{s^{\prime }}-\alpha }\left \Vert f\right \Vert
_{L_{p}\left( w^{p},B\left( x_{0},2r\right) \right) }\dint
\limits_{2r}^{\infty }\frac{dt}{t^{\frac{n}{s^{\prime }}-\alpha +1}} \\
&\leq &C\left( w^{q}\left( B\left( x_{0},r\right) \right) \right) ^{\frac{1}{%
q}}\left \Vert w^{-1}\right \Vert _{L_{\frac{s^{\prime }p}{\left(
p-s^{\prime }\right) }}\left( B\left( x_{0},r\right) \right) }\dint
\limits_{2r}^{\infty }\left \Vert f\right \Vert _{L_{p}\left( w^{p},B\left(
x_{0},t\right) \right) }\frac{dt}{t^{\frac{n}{s^{\prime }}-\alpha +1}} \\
&\leq &C\left( w^{q}\left( B\left( x_{0},r\right) \right) \right) ^{\frac{1}{%
q}}\dint \limits_{2r}^{\infty }\left \Vert f\right \Vert _{L_{p}\left(
w^{p},B\left( x_{0},t\right) \right) }\left \Vert w^{-1}\right \Vert _{L_{%
\frac{s^{\prime }p}{\left( p-s^{\prime }\right) }}\left( B\left(
x_{0},t\right) \right) }\frac{dt}{t^{\frac{n}{s^{\prime }}-\alpha +1}}.
\end{eqnarray*}%
Since $w\left( x\right) ^{s^{\prime }}\in A\left( \frac{p}{s^{\prime }},%
\frac{q}{s^{\prime }}\right) $, by (\ref{16}), we get%
\begin{equation}
\left( w^{q}\left( B\left( x_{0},t\right) \right) \right) ^{\frac{1}{q}%
}\left \Vert w^{-1}\right \Vert _{L_{\frac{s^{\prime }p}{\left( p-s^{\prime
}\right) }}\left( B\left( x_{0},t\right) \right) }\leq Ct^{\frac{n}{%
s^{\prime }}-\alpha }  \label{50}
\end{equation}%
holds for all $t>0$. Thus,%
\begin{eqnarray*}
\left \Vert T_{\Omega ,\alpha }^{A}f_{1}\right \Vert _{L_{q}\left(
w^{q},B\left( x_{0},r\right) \right) } &\leq &C\dsum \limits_{\left \vert
\gamma \right \vert =m-1}\left \Vert D^{\gamma }A\right \Vert _{\ast }\left(
w^{q}\left( B(x_{0},r)\right) \right) ^{\frac{1}{q}} \\
&&\times \dint \limits_{2r}^{\infty }\left( 1+\ln \frac{t}{r}\right) \left
\Vert f\right \Vert _{L_{p}\left( w^{p},B(x_{0},t)\right) }\left(
w^{q}\left( B(x_{0},t)\right) \right) ^{-\frac{1}{q}}\frac{1}{t}dt.
\end{eqnarray*}

Let $\Delta _{i}=B\left( x_{0},2^{j+1}r\right) \diagdown B\left(
x_{0},2^{j}r\right) $ and $x\in B(x_{0},r)$. By Lemma \ref{lemma5}, we get 
\begin{eqnarray*}
\left \vert T_{\Omega ,\alpha }^{A}f_{2}\left( x\right) \right \vert &\leq
&\left \vert \int \limits_{\left( B\left( x_{0},2r\right) \right) ^{C}}\frac{%
\Omega (x-y)}{|x-y|^{n-\alpha +m-1}}R_{m}\left( A;x,y\right) f(y)dy\right
\vert \\
&\leq &\dsum \limits_{j=1}^{\infty }\dint \limits_{\Delta _{i}}\frac{\left
\vert \Omega (x-y)f(y)\right \vert }{|x-y|^{n-\alpha }}\left( j\dsum
\limits_{\left \vert \gamma \right \vert =m-1}\left \Vert D^{\gamma }A\right
\Vert _{\ast }+\dsum \limits_{\left \vert \gamma \right \vert =m-1}\left
\vert D^{\gamma }A\left( y\right) -\left( D^{\gamma }A\right) _{B\left(
x_{0},r\right) }\right \vert \right) dy \\
&\leq &C\dsum \limits_{\left \vert \gamma \right \vert =m-1}\left \Vert
D^{\gamma }A\right \Vert _{\ast }\dsum \limits_{j=1}^{\infty }j\dint
\limits_{\Delta _{i}}\frac{\left \vert \Omega (x-y)f(y)\right \vert }{%
|x-y|^{n-\alpha }}dy \\
&&+C\dsum \limits_{\left \vert \gamma \right \vert =m-1}\dsum
\limits_{j=1}^{\infty }\dint \limits_{\Delta _{i}}\frac{\left \vert \Omega
(x-y)f(y)\right \vert }{|x-y|^{n-\alpha }}\left \vert D^{\gamma }A\left(
y\right) -\left( D^{\gamma }A\right) _{B\left( x_{0},r\right) }\right \vert
dy \\
&=&I_{1}+I_{2}.
\end{eqnarray*}%
By the H\"{o}lder's inequality, we have%
\begin{equation}
\dint \limits_{\Delta _{i}}\frac{\left \vert \Omega (x-y)f(y)\right \vert }{%
|x-y|^{n-\alpha }}dy\leq \left( \dint \limits_{\Delta _{i}}\left \vert
\Omega (x-y)\right \vert ^{s}dy\right) ^{\frac{1}{s}}\left( \dint
\limits_{\Delta _{i}}\frac{\left \vert f(y)\right \vert ^{s^{\prime }}}{%
|x-y|^{\left( n-\alpha \right) s^{\prime }}}dy\right) ^{\frac{1}{s^{\prime }}%
}.  \label{54}
\end{equation}

When $x\in B\left( x_{0},s\right) $ and $y\in \Delta _{i}$, then by a direct
calculation, we can see that $2^{j-1}r\leq \left \vert y-x\right \vert
<2^{j+1}r$. Hence,%
\begin{equation}
\left( \dint \limits_{\Delta _{i}}\left \vert \Omega (x-y)\right \vert
^{s}dy\right) ^{\frac{1}{s}}\leq C\left \Vert \Omega \right \Vert
_{L_{s}(S^{n-1})}\left \vert B\left( x_{0},2^{j+1}r\right) \right \vert ^{%
\frac{1}{s}}.  \label{55}
\end{equation}

It is clear that $x\in B\left( x_{0},r\right) $, $y\in B\left(
x_{0},2r\right) ^{C}$ implies \ $\frac{1}{2}\left \vert x_{0}-y\right \vert
\leq \left \vert x-y\right \vert \leq \frac{3}{2}\left \vert
x_{0}-y\right
\vert $. Consequently,%
\begin{equation*}
\left( \dint \limits_{\Delta _{i}}\frac{\left \vert f(y)\right \vert
^{s^{\prime }}}{|x-y|^{\left( n-\alpha \right) s^{\prime }}}dy\right) ^{%
\frac{1}{s^{\prime }}}\leq \frac{1}{\left \vert B\left(
x_{0},2^{j+1}r\right) \right \vert ^{1-\frac{\alpha }{n}}}\left( \dint
\limits_{B\left( x_{0},2^{j+1}r\right) }\left \vert f(y)\right \vert
^{s^{\prime }}dy\right) ^{\frac{1}{s^{\prime }}}.
\end{equation*}%
Then%
\begin{equation*}
I_{1}\leq C\dsum \limits_{\left \vert \gamma \right \vert =m-1}\left \Vert
D^{\gamma }A\right \Vert _{\ast }\dsum \limits_{j=1}^{\infty }j\left(
2^{j+1}r\right) ^{\alpha -\frac{n}{s^{\prime }}}\left( \dint
\limits_{B\left( x_{0},2^{j+1}r\right) }\left \vert f(y)\right \vert
^{s^{\prime }}dy\right) ^{\frac{1}{s^{\prime }}}.
\end{equation*}%
Since $s^{\prime }<p$, it follows from the H\"{o}lder's inequality that%
\begin{equation*}
\left( \dint \limits_{B\left( x_{0},2^{j+1}r\right) }\left \vert f(y)\right
\vert ^{s^{\prime }}dy\right) ^{\frac{1}{s^{\prime }}}\leq C\left \Vert
f\right \Vert _{L_{p}\left( w^{p},B\left( x_{0},2^{j+1}r\right) \right)
}\left \Vert w^{-1}\right \Vert _{L_{\frac{s^{\prime }p}{\left( p-s^{\prime
}\right) }}\left( B\left( x_{0},2^{j+1}r\right) \right) }.
\end{equation*}%
Then%
\begin{eqnarray*}
&&\dsum \limits_{j=1}^{\infty }j\left( 2^{j+1}r\right) ^{\alpha -\frac{n}{%
s^{\prime }}}\left( \dint \limits_{B\left( x_{0},2^{j+1}r\right) }\left
\vert f(y)\right \vert ^{s^{\prime }}dy\right) ^{\frac{1}{s^{\prime }}} \\
&\leq &C\dsum \limits_{j=1}^{\infty }\left( 1+\ln \frac{2^{j+1}r}{r}\right)
\left( 2^{j+1}r\right) ^{\alpha -\frac{n}{s^{\prime }}}\left \Vert f\right
\Vert _{L_{p}\left( w^{p},B\left( x_{0},2^{j+1}r\right) \right) }\left \Vert
w^{-1}\right \Vert _{L_{\frac{s^{\prime }p}{\left( p-s^{\prime }\right) }%
}\left( B\left( x_{0},2^{j+1}r\right) \right) } \\
&\leq &C\dsum \limits_{j=1}^{\infty }\dint
\limits_{2^{j+1}r}^{2^{j+2}r}\left( 1+\ln \frac{t}{r}\right) \left \Vert
f\right \Vert _{L_{p}\left( w^{p},B(x_{0},t)\right) }\left \Vert
w^{-1}\right \Vert _{L_{\frac{s^{\prime }p}{\left( p-s^{\prime }\right) }%
}\left( B\left( x_{0},t\right) \right) }\frac{dt}{t^{\frac{n}{s^{\prime }}%
-\alpha +1}} \\
&\leq &C\dint \limits_{2r}^{\infty }\left( 1+\ln \frac{t}{r}\right) \left
\Vert f\right \Vert _{L_{p}\left( w^{p},B(x_{0},t)\right) }\left \Vert
w^{-1}\right \Vert _{L_{\frac{s^{\prime }p}{\left( p-s^{\prime }\right) }%
}\left( B\left( x_{0},t\right) \right) }\frac{dt}{t^{\frac{n}{s^{\prime }}%
-\alpha +1}}.
\end{eqnarray*}%
By (\ref{50}), we know%
\begin{equation}
\left \Vert w^{-1}\right \Vert _{L_{\frac{s^{\prime }p}{\left( p-s^{\prime
}\right) }}\left( B\left( x_{0},t\right) \right) }\leq Ct^{\frac{n}{%
s^{\prime }}-\alpha }\left( w^{q}\left( B(x_{0},t)\right) \right) ^{-\frac{1%
}{q}}.  \label{59}
\end{equation}%
Then%
\begin{equation*}
I_{1}\leq C\dsum \limits_{\left \vert \gamma \right \vert =m-1}\left \Vert
D^{\gamma }A\right \Vert _{\ast }\dint \limits_{2r}^{\infty }\left( 1+\ln 
\frac{t}{r}\right) \left \Vert f\right \Vert _{L_{p}\left(
w^{p},B(x_{0},t)\right) }\left( w^{q}\left( B(x_{0},t)\right) \right) ^{-%
\frac{1}{q}}\frac{1}{t}dt.
\end{equation*}

On the other hand, by the H\"{o}lder's inequality, (\ref{54}) and (\ref{55})
we have%
\begin{eqnarray*}
&&\dint \limits_{\Delta _{i}}\frac{\left \vert \Omega (x-y)f(y)\right \vert 
}{|x-y|^{n-\alpha }}\left \vert D^{\gamma }A\left( y\right) -\left(
D^{\gamma }A\right) _{B\left( x_{0},r\right) }\right \vert dy \\
&\leq &\left( \dint \limits_{\Delta _{i}}\left \vert \Omega (x-y)\right
\vert ^{s}dy\right) ^{\frac{1}{s}}\left( \dint \limits_{\Delta _{i}}\frac{%
\left \vert D^{\gamma }A\left( y\right) -\left( D^{\gamma }A\right)
_{B\left( x_{0},r\right) }f(y)\right \vert ^{s^{\prime }}}{|x-y|^{\left(
n-\alpha \right) s^{\prime }}}dy\right) ^{\frac{1}{s^{\prime }}} \\
&\leq &C\dsum \limits_{j=1}^{\infty }\left( 2^{j+1}r\right) ^{\alpha -\frac{n%
}{s^{\prime }}}\left( \dint \limits_{B\left( x_{0},2^{j+1}r\right) }\left
\vert D^{\gamma }A\left( y\right) -\left( D^{\gamma }A\right) _{B\left(
x_{0},r\right) }\right \vert ^{s^{\prime }}\left \vert f(y)\right \vert
^{s^{\prime }}dy\right) ^{\frac{1}{s^{\prime }}}.
\end{eqnarray*}%
Applying the H\"{o}lder's inequality again, we get%
\begin{eqnarray*}
&&\left( \dint \limits_{B\left( x_{0},2^{j+1}r\right) }\left \vert D^{\gamma
}A\left( y\right) -\left( D^{\gamma }A\right) _{B\left( x_{0},r\right)
}\right \vert ^{s^{\prime }}\left \vert f(y)\right \vert ^{s^{\prime
}}dy\right) ^{\frac{1}{s^{\prime }}} \\
&\leq &C\left \Vert f\right \Vert _{L_{p}\left( w^{p},B\left(
x_{0},2^{j+1}r\right) \right) }\left \Vert \left( D^{\gamma }A\left(
y\right) -\left( D^{\gamma }A\right) _{B\left( x_{0},r\right) }\right)
w\left( \cdot \right) ^{-1}\right \Vert _{L_{\frac{s^{\prime }p}{\left(
p-s^{\prime }\right) }}\left( B\left( x_{0},2^{j+1}r\right) \right) }.
\end{eqnarray*}%
Consequently,%
\begin{eqnarray*}
I_{2} &\leq &C\dsum \limits_{\left \vert \gamma \right \vert =m-1}\left
\Vert D^{\gamma }A\right \Vert _{\ast }\dsum \limits_{j=1}^{\infty }\dint
\limits_{2^{j+1}r}^{2^{j+2}r}\left( 1+\ln \frac{2^{j+1}r}{r}\right) \left(
2^{j+1}r\right) ^{\alpha -\frac{n}{s^{\prime }}}\left \Vert f\right \Vert
_{L_{p}\left( w^{p},B(x_{0},t)\right) } \\
&&\times \left \Vert \left( D^{\gamma }A\left( y\right) -\left( D^{\gamma
}A\right) _{B\left( x_{0},r\right) }\right) w\left( \cdot \right)
^{-1}\right \Vert _{L_{\frac{s^{\prime }p}{\left( p-s^{\prime }\right) }%
}\left( B\left( x_{0},t\right) \right) }dt \\
&\leq &C\dsum \limits_{\left \vert \gamma \right \vert =m-1}\left \Vert
D^{\gamma }A\right \Vert _{\ast }\dint \limits_{2r}^{\infty }\left \Vert
f\right \Vert _{L_{p}\left( w^{p},B(x_{0},t)\right) } \\
&&\times \left \Vert \left( D^{\gamma }A\left( y\right) -\left( D^{\gamma
}A\right) _{B\left( x_{0},r\right) }\right) w\left( \cdot \right)
^{-1}\right \Vert _{L_{\frac{s^{\prime }p}{\left( p-s^{\prime }\right) }%
}\left( B\left( x_{0},t\right) \right) }\frac{dt}{t^{\frac{n}{s^{\prime }}%
-\alpha +1}}.
\end{eqnarray*}

By $w\left( x\right) ^{s^{\prime }}\in A\left( \frac{p}{s^{\prime }},\frac{q%
}{s^{\prime }}\right) $ and $\left( ii\right) $ of Lemma \ref{lemma0}, we
know $w(x)^{-\frac{s^{\prime }p}{\left( p-s^{\prime }\right) }}\in A_{1+%
\frac{s^{\prime }p}{\left( p-s^{\prime }\right) q}}$. Then it follows from (%
\ref{31}) and (\ref{59}) that%
\begin{eqnarray*}
&&\left \Vert \left( D^{\gamma }A\left( y\right) -\left( D^{\gamma }A\right)
_{B\left( x_{0},r\right) }\right) w\left( \cdot \right) ^{-1}\right \Vert
_{L_{\frac{s^{\prime }p}{\left( p-s^{\prime }\right) }}\left( B\left(
x_{0},t\right) \right) }\frac{dt}{t^{\frac{n}{s^{\prime }}-\alpha +1}} \\
&\leq &\left( \dint \limits_{B\left( x_{0},r\right) }\left \vert \left(
D^{\gamma }A\left( y\right) -\left( D^{\gamma }A\right) _{B\left(
x_{0},r\right) }\right) \right \vert ^{\frac{s^{\prime }p}{\left(
p-s^{\prime }\right) }}w^{-\frac{s^{\prime }p}{\left( p-s^{\prime }\right) }%
}\left( y\right) dy\right) ^{\frac{\left( p-s^{\prime }\right) }{s^{\prime }p%
}} \\
&\leq &C\left \Vert D^{\gamma }A\right \Vert _{\ast }\left( 1+\ln \frac{t}{r}%
\right) \left( w^{-\frac{s^{\prime }p}{\left( p-s^{\prime }\right) }}\left(
B\left( x_{0},t\right) \right) \right) ^{\frac{\left( p-s^{\prime }\right) }{%
s^{\prime }p}} \\
&=&C\left \Vert D^{\gamma }A\right \Vert _{\ast }\left( 1+\ln \frac{t}{r}%
\right) \left \Vert w^{-1}\right \Vert _{L_{\frac{s^{\prime }p}{\left(
p-s^{\prime }\right) }}\left( B\left( x_{0},t\right) \right) } \\
&\leq &C\left \Vert D^{\gamma }A\right \Vert _{\ast }\left( 1+\ln \frac{t}{r}%
\right) t^{\frac{n}{s^{\prime }}-\alpha }\left( w^{q}\left(
B(x_{0},t)\right) \right) ^{-\frac{1}{q}}.
\end{eqnarray*}%
Thus,%
\begin{equation*}
I_{2}\leq C\dsum \limits_{\left \vert \gamma \right \vert =m-1}\left \Vert
D^{\gamma }A\right \Vert _{\ast }\dint \limits_{2r}^{\infty }\left( 1+\ln 
\frac{t}{r}\right) \left \Vert f\right \Vert _{L_{p}\left(
w^{p},B(x_{0},t)\right) }\left( w^{q}\left( B(x_{0},t)\right) \right) ^{-%
\frac{1}{q}}\frac{1}{t}dt.
\end{equation*}

Combining the estimates of $I_{1}$ and $I_{2}$, we get 
\begin{equation*}
\left \vert T_{\Omega ,\alpha }^{A}f_{2}\left( x\right) \right \vert \leq
C\dsum \limits_{\left \vert \gamma \right \vert =m-1}\left \Vert D^{\gamma
}A\right \Vert _{\ast }\dint \limits_{2r}^{\infty }\left( 1+\ln \frac{t}{r}%
\right) \left \Vert f\right \Vert _{L_{p}\left( w^{p},B(x_{0},t)\right)
}\left( w^{q}\left( B(x_{0},t)\right) \right) ^{-\frac{1}{q}}\frac{1}{t}dt.
\end{equation*}%
Then we get 
\begin{eqnarray*}
\left \Vert T_{\Omega ,\alpha }^{A}f_{2}\right \Vert _{L_{q}\left(
w^{q},B\left( x_{0},r\right) \right) } &\leq &C\dsum \limits_{\left \vert
\gamma \right \vert =m-1}\left \Vert D^{\gamma }A\right \Vert _{\ast }\left(
w^{q}\left( B(x_{0},r)\right) \right) ^{\frac{1}{q}} \\
&&\times \dint \limits_{2r}^{\infty }\left( 1+\ln \frac{t}{r}\right) \left
\Vert f\right \Vert _{L_{p}\left( w^{p},B(x_{0},t)\right) }\left(
w^{q}\left( B(x_{0},t)\right) \right) ^{-\frac{1}{q}}\frac{1}{t}dt.
\end{eqnarray*}%
This completes the proof of Theorem \ref{teo2}.

\textbf{Proof of Theorem \ref{teo3}.}

We consider (\ref{5}) firstly. Since $f\in M_{p,\varphi _{1}}\left( w^{p},{%
\mathbb{R}^{n}}\right) $, by (\ref{11}) and the fact that $\left \Vert
f\right \Vert _{L_{p}\left( w^{p},B\left( x_{0},t\right) \right) }$ is a
non-decreasing function of $t$, we get%
\begin{eqnarray*}
&&\frac{\left \Vert f\right \Vert _{L_{p}\left( w^{p},B\left( x_{0},t\right)
\right) }}{\limfunc{essinf}\limits_{0<t<\tau <\infty }\varphi
_{1}(x_{0},\tau )\left( w^{p}\left( B(x_{0},\tau )\right) \right) ^{\frac{1}{%
p}}} \\
&\leq &\limfunc{esssup}\limits_{0<t<\tau <\infty }\frac{\left \Vert f\right
\Vert _{L_{p}\left( w^{p},B\left( x_{0},t\right) \right) }}{\varphi
_{1}(x_{0},\tau )\left( w^{p}\left( B(x_{0},\tau )\right) \right) ^{\frac{1}{%
p}}} \\
&\leq &\limfunc{esssup}\limits_{0<\tau <\infty }\frac{\left \Vert f\right
\Vert _{L_{p}\left( w^{p},B\left( x_{0},\tau \right) \right) }}{\varphi
_{1}(x_{0},\tau )\left( w^{p}\left( B(x_{0},\tau )\right) \right) ^{\frac{1}{%
p}}} \\
&\leq &\left \Vert f\right \Vert _{M_{p,\varphi _{1}}\left( w^{p},{\mathbb{R}%
^{n}}\right) }.
\end{eqnarray*}

For $s^{\prime }<p<\infty $, since $(\varphi _{1},\varphi _{2})$ satisfies (%
\ref{4}), we have%
\begin{eqnarray*}
&&\int \limits_{r}^{\infty }\left( 1+\ln \frac{t}{r}\right) \left \Vert
f\right \Vert _{L_{p}\left( w^{p},B\left( x_{0},t\right) \right) }\left(
w^{q}\left( B(x_{0},t)\right) \right) ^{-\frac{1}{q}}\frac{dt}{t} \\
&\leq &\int \limits_{r}^{\infty }\left( 1+\ln \frac{t}{r}\right) \frac{\left
\Vert f\right \Vert _{L_{p}\left( w^{p},B\left( x_{0},t\right) \right) }}{%
\limfunc{essinf}\limits_{t<\tau <\infty }\varphi _{1}(x_{0},\tau )\left(
w^{p}\left( B(x_{0},\tau )\right) \right) ^{\frac{1}{p}}}\frac{\limfunc{%
essinf}\limits_{t<\tau <\infty }\varphi _{1}(x_{0},\tau )\left( w^{p}\left(
B(x_{0},\tau )\right) \right) ^{\frac{1}{p}}}{\left( w^{q}\left(
B(x_{0},t)\right) \right) ^{\frac{1}{q}}}\frac{dt}{t} \\
&\leq &C\left \Vert f\right \Vert _{M_{p,\varphi _{1}}\left( w^{p},{\mathbb{R%
}^{n}}\right) }\int \limits_{r}^{\infty }\left( 1+\ln \frac{t}{r}\right) 
\frac{\limfunc{essinf}\limits_{t<\tau <\infty }\varphi _{1}(x_{0},\tau
)\left( w^{p}\left( B(x_{0},\tau )\right) \right) ^{\frac{1}{p}}}{\left(
w^{q}\left( B(x_{0},t)\right) \right) ^{\frac{1}{q}}}\frac{dt}{t}
\end{eqnarray*}%
\begin{equation}
\leq C\left \Vert f\right \Vert _{M_{p,\varphi _{1}}\left( w^{p},{\mathbb{R}%
^{n}}\right) }\varphi _{2}(x_{0},r).  \label{0}
\end{equation}%
Then by (\ref{3}) and (\ref{0}), we get%
\begin{eqnarray*}
\left \Vert T_{\Omega ,\alpha }^{A}f\right \Vert _{M_{q,\varphi _{2}}\left(
w^{q},{\mathbb{R}^{n}}\right) } &=&\sup_{x_{0}\in {\mathbb{R}^{n},}%
r>0}\varphi _{2}\left( x_{0},r\right) ^{-1}\left( w^{q}\left(
B(x_{0},r)\right) \right) ^{-\frac{1}{q}}\left \Vert T_{\Omega ,\alpha
}^{A}f\right \Vert _{L_{q}\left( w^{q},B(x_{0},r)\right) } \\
&\leq &C\dsum \limits_{\left \vert \gamma \right \vert =m-1}\left \Vert
D^{\gamma }A\right \Vert _{\ast }\sup_{x_{0}\in {\mathbb{R}^{n},}r>0}\varphi
_{2}\left( x_{0},r\right) ^{-1} \\
&&\times \dint \limits_{r}^{\infty }\left( 1+\ln \frac{t}{r}\right) \left
\Vert f\right \Vert _{L_{p}\left( w^{p},B(x_{0},t)\right) }\left(
w^{q}\left( B(x_{0},t)\right) \right) ^{-\frac{1}{q}}\frac{1}{t}dt \\
&\leq &C\dsum \limits_{\left \vert \gamma \right \vert =m-1}\left \Vert
D^{\gamma }A\right \Vert _{\ast }\left \Vert f\right \Vert _{M_{p,\varphi
_{1}}\left( w^{p},{\mathbb{R}^{n}}\right) }.
\end{eqnarray*}%
Hence, we have completed the proof of (\ref{5}).

We are now in a place of proving (\ref{6}) in Theorem \ref{teo3}.

Set%
\begin{equation*}
\widetilde{T}_{\left \vert \Omega \right \vert ,\alpha }^{A}\left( \left
\vert f\right \vert \right) (x)=\int \limits_{{\mathbb{R}^{n}}}\frac{\left
\vert \Omega (x-y)\right \vert }{|x-y|^{n-\alpha +m-1}}\left \vert
R_{m}\left( A;x,y\right) \right \vert \left \vert f(y)\right \vert dy\qquad
0<\alpha <n,
\end{equation*}%
where $\Omega \in L_{s}(S^{n-1})\left( s>1\right) $ is homogeneous of degree
zero in ${\mathbb{R}^{n}}$. It is easy to see that, for $\widetilde{T}%
_{\left \vert \Omega \right \vert ,\alpha }^{A}$, the conclusions of
inequality (\ref{5}) also hold. On the other hand, for any $r>0$, we have 
\begin{eqnarray*}
\widetilde{T}_{\left \vert \Omega \right \vert ,\alpha }^{A}\left( \left
\vert f\right \vert \right) (x) &\geq &\dint \limits_{|x-y|<r}\frac{\left
\vert \Omega (x-y)\right \vert }{|x-y|^{n-\alpha +m-1}}\left \vert
R_{m}\left( A;x,y\right) \right \vert \left \vert f(y)\right \vert dy \\
&\geq &\frac{1}{r^{n-\alpha +m-1}}\dint \limits_{|x-y|<r}\left \vert \Omega
(x-y)\right \vert \left \vert R_{m}\left( A;x,y\right) \right \vert \left
\vert f(y)\right \vert dy.
\end{eqnarray*}%
Taking the supremum for $r>0$ on the inequality above, we get%
\begin{equation}
\widetilde{T}_{\left \vert \Omega \right \vert ,\alpha }^{A}\left( \left
\vert f\right \vert \right) (x)\geq M_{\Omega ,\alpha }^{A}f(x)\qquad \text{%
for }x\in {\mathbb{R}^{n}}.  \label{100}
\end{equation}%
Thus, we can immediately obtain (\ref{6}) from (\ref{100}) and (\ref{5}),
which completes the proof.

\end{document}